\documentclass[10pt]{elsarticle}
\usepackage{amsmath}
\usepackage{amssymb}
\usepackage{pifont}
\usepackage{amsfonts}
\usepackage{graphicx}
\usepackage{setspace}
\usepackage[left=2cm,top=1.75cm,right=2cm]{geometry}

\newtheorem{lemma}{Lemma}
\newtheorem{remark}{Remark}
\newtheorem{theorem}{Theorem}

\numberwithin{equation}{section}

\begin{document}
\begin{frontmatter}
\title{On approximation properties of Baskakov-Schurer-Sz\'{a}sz operators}
\author[label1,label2,label*]{Vishnu Narayan Mishra}
\ead{vishnunarayanmishra@gmail.com,
vishnu\_narayanmishra@yahoo.co.in}
\address[label1]{Department of Applied Mathematics \& Humanities,
Sardar Vallabhbhai National Institute of Technology, Ichchhanath Mahadev Dumas Road, Surat -395 007 (Gujarat), India}
\address[label2]{L. 1627 Awadh Puri Colony Beniganj, Phase - III, Opposite - Industrial Training Institute, Ayodhya Main Road, Faizabad-224 001,(Uttar Pradesh), India }
\fntext[label*]{Corresponding author}
\author[label1]{Preeti Sharma}
\ead{preeti.iitan@gmail.com}

\begin{abstract}
In this paper, we are dealing with a new type of Baskakov-Schurer-Sz\'{a}sz operators (\ref{eq1}).
Approximation properties of this operators are explored: the rate of convergence in terms of the usual moduli of smoothness is given, the convergence in certain weighted spaces is investigated. We study $q$-analogues of Baskakov-Schurer-Sz\'{a}sz operators and it's Stancu generalization. In the last section, we give better error estimations for the operators (\ref{eq2}) using King type approach and obtained weighted statistical approximation properties for operator (\ref{qb}).
\end{abstract}

\end{frontmatter}

\section{Introduction}
This work is at confluence of two mathematical research areas namely linear approximation processes and statistical
convergence. To reveal what novelties this paper brings, we briefly present both a background and some historical comments. The main problem of approximation theory consists in finding for a complicated function a close-by simple function.\\
\indent At 70 years old, Karl Wilhelm Theodor Weierstrass (1815--1897) proved the density of the algebraic polynomials in the space $C([a,b])$ and of the trigonometric polynomial in $\tilde{C}([a,b]).$
Weierstrass’s approximation theorem stating that every continuous function on a bounded interval can be approximated to arbitrary accuracy by polynomials is such an important example for this process and has been played
the significant role in the development of analysis. By using probability theory Bernstein \cite{B} proved the
Weierstrass’s theorem and defined approximate polynomials known as Bernstein polynomials in the literature.
 In 1987, Lupa\c{s} \cite{L} introduced the first $q$-analogue
of Bernstein operator and investigated its approximating and shape-preserving properties.
Another $q$-generalization of the classical Bernstein polynomials is due to Phillips
\cite{PH}. After that many generalizations of well-known positive linear operators, based
on $q$-integers were introduced and studied by several authors. In \cite{OA} Agratini introduced
a new class of $q$-Bernstein-type operators which fix certain polynomials and
studied the limit of iterates of Lupa\c{s} $q$-analogue of the Bernstein operators.
In 1974, Khan \cite{HH} studied approximation of functions in various classes using different
types of operators.


In 1950, Sz\'{a}sz defined and studied the approximation properties of the following operators
\begin{equation}
S_n(f;x)=e^{-n x}\sum_{j=0}^{\infty} \frac{(n x)^j}{j!}f\bigg(\frac{j}{n}\bigg),
\end{equation}
whenever $f$ satisfies exponential-type growth condition \cite{OS}. In 1962, Schurer \cite{FS} introduced and studied the approximation properties of linear positive operators. An extension in $q$-Calculus of Sz\'{a}sz-Mirakyan operators was constructed by Aral \cite{AA} who formulated also a Voronovskaya theorem related to $q$- derivatives for these operators. In 1976, May \cite{M} showed that the Baskakov operators can reduce to the Sz\'{a}sz-Mirakyan operators.  After that several other researchers have studied in this direction and obtained different approximation properties of many operators \cite{KJ}-\cite{MA}.

Motivated by the operators due to Y\"{u}ksel \cite{IY}, in the year 2013, we considered the following
operators.
\newline
\indent Let $p, k, n \in \mathbb{N}$ and $f$ be a real valued continuous function on the interval $[0,\infty)$. We define the Baskakov-Schurer-Sz\'{a}sz type linear positive operators as 

\begin{equation}\label{eq1}
\mathcal{L}_{n,p}(f,x)= (n+p) \sum_{k=0}^{\infty} b^k_{n,p}(x)\int_{0}^{\infty} f(t) s^k_{n,p}(t) d t \\
\end{equation}


where  $\displaystyle b^{k}_{n,p}(x)=\binom{n+p+k-1}{k}\frac{x^k}{(1+x)^{n+p+k}}$\,\,
and\,\, $\displaystyle s^k_{n,p}(t)= e^{-(n+p)t}\frac{{((n+p)t)}^k}{k!}.$\\
\newline
Now we give an auxiliary lemma for the Korovkin test functions.
\begin{lemma}
Let $e_m(t) =t^m$, $m = 0,1,2.$ we have
\begin{itemize}
\item[(i)] $\mathcal{L}_{n,p}(e_0;x)=1,$
\item[(ii)] $\mathcal{L}_{n,p}(e_1;x)=x+\frac{1}{(n+p)},$
\item[(iii)] $\mathcal{L}_{n,p}(e_2;x)=\big(1+\frac{1}{n+p}\big)x^2+\frac{4}{(n+p)}x+\frac{2}{(n+p)^2}.$
\end{itemize}
\end{lemma}

\begin{lemma}\label{l1}
Let $p\in \mathbb{N}$. Then following hold:\\
$\alpha_n(x):= \mathcal{L}_{n,p}((t-x);x)=\frac{1}{n+p}.$\\
$\delta_n(x):=\mathcal{L}_{n,p}((t-x)^2;x)=\frac{1}{n+p} x^2+\frac{2}{n+p} x+ \frac{2}{(n+p)^2}.$
\end{lemma}

\section{Some auxiliary results}
Let the space $C_B[0,\infty)$ of all continuous and bounded
functions $f$ on $[0,\infty)$, be endowed with the norm $\|f\|=sup\{ \mid f(x)\mid: x\in[0,\infty)\}.$ Further let us consider the Peetre's K-functional which is defined by
 \begin{equation} K_2(f,\delta)= \inf_{g\in W^2}\{\|f-g\|+\delta \|g''\|\},\end{equation}
where $\delta >0$ and
$W^2_{\infty}=\{g\in C_B[0,\infty):g', g'' \in C_B[0,\infty)\}.$ By the method
as given (\cite{DL} p.177, Theorem 2.4), there exists an absolute constant $C>0$ such that
\begin{equation}\label{1}
K_2(f,\delta) \leq C \omega_2(f,\delta),\end{equation}
 where
 
\begin{equation}\label{2}
\omega_2(f,{\delta})=\sup_{0<h \leq \delta}\sup_{x\in[0,\infty)}\mid f(x+2h)-2f(x+h)+f(x)\mid \end{equation}
is the second order modulus of smoothness of $f\in C_B[0,\infty).$
Also we set

\begin{equation}\label{3}
\omega(f,{\delta})= \sup_{0<h \leq
\delta}\sup_{x\in[0,\infty)}\mid f(x+h)-f(x)\mid.
\end{equation}
We denote the usual modulus of continuity of $f\in C_B[0,\infty)$.

\begin{theorem}\label{the1}
Let $f \in C_B[0,\infty)$, then for all $x\in
[0,\infty)$, there exists an absolute constant $ C >0 $ such that
\begin{equation}|\mathcal{L}_{n,p}(f,x)-f(x)| \leq C \omega_2\left(f,\sqrt{\delta_n(x) +(\alpha_n(x))^2}\right)+ 
\omega(f,\alpha_n(x)).
\end{equation}
\end{theorem}
\textbf{Proof.}
Let $g\in W_\infty^2$ and $x,t\in [0,\infty).$ By Taylor's expansion, we have
\begin{equation} \label{ac1}
g(t)=g(x)+g'(x)(t-x)+\int_x^t (t-u)g''(u)du.
\end{equation}
Define \begin{equation}
\mathcal{\widetilde{L}}_{n,p}(f,x) = \mathcal{L}_{n,p}(f,x) + f(x)-f\left(x +\frac{1}{n+p}\right).
\end{equation}
Now, we have $\mathcal{\widetilde{L}}_{n,p}(t-x,x)=0,$ $t\in[0,\infty).$\\
Applying $\mathcal{\widetilde{L}}_{n,p}$ on both sides of (\ref{ac1}), we get
\begin{eqnarray*}
\mathcal{\widetilde{L}}_{n,p}(g,x)-g(x) & = & g'(x)\mathcal{\tilde{L}}_{n,p}((t-x),x)+\mathcal{\tilde{L}}_{n,p} \left(\int_x^t(t-u)g''(u)du,x\right)\\
&=& \mathcal{L}_{n,p}\left(\int_x^t(t-u)g''(u) du, x\right)+\int_x^{x +\frac{1}{n+p}} \left(x +\frac{1}{n+p}-u\right)g''(u)du,
\end{eqnarray*}
on the other hand $\displaystyle \bigg| \int_x^t (x-u)g''(u)du\bigg| \leq \|g''\|(t-x)^2$ and
\begin{eqnarray*}
\int_x^{x +\frac{1}{n+p}} \left(x +\frac{1}{n+p}-u\right)g''(u)du \leq \left(x +\frac{1}{n+p}-x\right)^2\|g''\|&=& \|g''\|\left(\mathcal{L}_{n,p}(t-x,x)\right)^2\\
&=&\|g''\|\left(\alpha_n(x)\right)^2.
\end{eqnarray*}
Thus, one can do this
\begin{eqnarray*}
\bigg|\mathcal{\tilde{L}}_{n,p}(g,x)-g(x)\bigg| &\leq & \bigg|\mathcal{L}_{n,p} \left(\int_x^t(t-u)g''(u)du,x\right)\bigg| +\bigg|\int_x^{x +\frac{1}{n+p}}\left(x +\frac{1}{n+p} -u\right)g''(u)du\bigg|\\
&\leq& \|g''\|\mathcal{L}_{n,p} \left((t-x)^2,x\right)+\|g''\|\left(\alpha_n(x)\right)^2\\
&\leq& \|g''\| \bigg[\delta_n(x) +\left(\alpha_n(x)\right)^2\bigg].
\end{eqnarray*}
We observe that,
\begin{eqnarray*}
\bigg|\mathcal{\widetilde{L}}_{n,p}(f,x)-f(x)\bigg| &\leq & \bigg|\mathcal{\tilde{L}}_{n,p}(f-g,x)-(f-g)(x)\bigg|\\
&+&  \bigg|\mathcal{\tilde{L}}_{n,p} (g,x)-g(x)\bigg|+\bigg|f(x)-f\left(x +\frac{1}{n+p}\right)\bigg|\\
 &\leq & \|f-g\|+\|g''\| \bigg[\delta_n(x)+\left(\alpha_n(x)\right)^2\bigg]+\omega(f,\alpha_n(x)).
\end{eqnarray*}
Now, taking infimum on the right-hand side over all $g\in C_B^2[0,\infty)$  and from (\ref{1}), we get
\begin{eqnarray*}
\bigg|\mathcal{L}_{n,p}(f,x)-f(x)\bigg| &\leq & K_2\left(f,\delta_{n}(x)+\left(\alpha_{n}(x)\right)^2\right)+ \omega(f,\alpha_n(x))\\
&\leq& C \omega_2\left(f,\sqrt{\delta_{n}(x)+\left(\alpha_{n}(x)\right)^2}\right)+ \omega(f,\alpha_{n}(x)),
\end{eqnarray*}
which proves the theorem.\\

\section{Weighted approximation}
In this section, we obtain the Korovkin type weighted approximation by the operators defined in (\ref{eq1}). The weighted Korovkin-type theorems were proved by Gadzhiev \cite{AD}.
A real function $\rho = 1+x^2$ is called a weight function if it is continuous on $%
\mathbb{R}$ and $\lim\limits_{\mid x\mid \rightarrow \infty }\rho (x)=\infty
,~\rho (x)\geq 1$ for all $x\in \mathbb{R}$.

Let $B_{\rho }(\mathbb{R})$  denote the weighted space of real-valued
functions $f$ defined on $\mathbb{R}$ with the property $\mid f(x)\mid \leq
M_{f}~\rho (x)$ for all $x\in \mathbb{R}$, where $M_{f}$ is a constant
depending on the function $f$. We also consider the weighted subspace $%
C_{\rho }(\mathbb{R})$ of $B_{\rho }(\mathbb{R})$ given by $C_{\rho }(%
\mathbb{R})=\{f\in B_{\rho }(\mathbb{R}){:}$ $f$ is continuous on $\mathbb{R} $\} and $C_{\rho}^{*}[0,\infty)$ denotes the subspace of all functions
 $f\in C_{\rho}[0,\infty)$ for which $\lim\limits_{|x|\rightarrow\infty} \frac{ f(x)}{\rho(x)}$ exists finitely.

\begin{theorem}(See \cite{AD} and \cite{AD1})
\begin{itemize}
\item[(i)] There exists a sequence of linear positive operators $A_n(C_{\rho}\rightarrow B_{\rho})$ such that
\begin{equation}\label{wa1}
\lim_{n\rightarrow \infty}\|A_n(\phi^\nu)- \phi^\nu\|_{\rho}=0,~~\nu=0,1,2
\end{equation}
and a function $f^{*}\in C_{\rho} \backslash C^{*}_{\rho}$  with $\lim\limits_{n\rightarrow \infty} \| A_n(f^{*})- f^{*}\|_{\rho}\geq 1.$
\item[(ii)] If a sequence of linear positive operators $A_n(C_{\rho}\rightarrow B_{\rho})$ satisfies conditions (\ref{wa1}) then
\begin{equation}
\lim_{n\rightarrow \infty}\|A_n(f)- f\|_{\rho}=0, \text{  for every } f\in C^{*}_{\rho}.
\end{equation}
\end{itemize}
\end{theorem}
Throughout this paper we take the growth condition as $\rho(x) = 1 + x^2$ and $\rho_{\gamma}(x) = 1 + x^{2+\gamma},~ x\in [0,\infty), \gamma > 0.$
Now we are ready to prove our next result as follows:

\begin{theorem}
For each $f \in C_{\rho}^{*}[0,\infty)$, we have
\begin{equation*}\lim _{n\rightarrow \infty} \| \mathcal{L}_{n,p}(f)-f \|_{\rho} = 0.\end{equation*}
\end{theorem}
\textbf{Proof.} Using the theorem in \cite{AD} we see that it is sufficient to verify the following three conditions
\begin{equation}\label{w1}
\lim_{n\rightarrow \infty} \|
\mathcal{L}_{n,p}(t^r,x)-x^r\|_{\rho}=0, \text{  } r=0,1,2.
\end{equation}
Since, $\mathcal{L}_{n,p}(1,x)=1$, the first condition of (\ref{w1}) is satisfied for $r=0$. Now,
\begin{eqnarray*}
\|\mathcal{L}_{n,p}(t,x)-x\|_{\rho}&=&\sup_{x\in [0,\infty)} \frac{\mid \mathcal{L}_{n,p}(t,x)-x \mid}{1+x^2}\\
&\leq&  \frac{1}{n+p} \sup_{x\in [0,\infty)}\frac{1}{1+x^2}\\
&=& o(1).
\end{eqnarray*}
Finally,
\begin{eqnarray*}
\|\mathcal{L}_{n,p}(t^2,x)-x^2\|_{\rho} & = & \sup_{x\in [0,\infty)} \frac{\mid \mathcal{L}_{n,p}(t^2,x)-x^2 \mid}{1+x^2}\\
&\leq&\left|\frac{1}{(n+p)}\right|\sup_{x\in [0,\infty)} \frac{x^2}{1+x^2}+
\left|\frac{4}{(n+p)}\right|\sup_{x\in [0,\infty)}\frac{x}{1+x^2}+\left| \frac{2}{(n+p)^2}\right| \\
&=& o(1).
\end{eqnarray*}
Thus, from Gadziv's Theorem in \cite{AD} we obtain the desired result of theorem.
\qed

We give the following theorem to approximate all functions in $C_{x^2}[0,\infty)$.
\begin{theorem}
For each $f\in C_{x^2}[0,\infty)$ and $\alpha
>0$, we have \\
$$\lim\limits_{n \rightarrow \infty} \sup_{x\in[0,\infty)} \frac{\mid\mathcal{L}_{n,p}(f,x)-f(x)
\mid}{(1+x^2)^{1+\alpha}}=0.$$
\end{theorem}
\textbf{Proof.} For any fixed $x_0>0$,
\begin{eqnarray*}
\sup_{x\in[0,\infty)}\frac{\mid \mathcal{L}_{n,p}(f,x)-f(x) \mid}{(1+x^2)^{1+\alpha}}&\leq & \sup_{x \leq x_0}\frac{\mid \mathcal{L}_{n,p}(f,x)-f(x) \mid}{(1+x^2)^{1+\alpha}} + \sup_{x \geq x_0}\frac{\mid \mathcal{L}_{n,p}(f,x)-f(x) \mid}{(1+x^2)^{1+\alpha}}\\
&\leq&  \|\mathcal{L}_{n,p}(f)-f\|_{C[0,x_0]} + \|f\|_{x^2}\sup_{x \geq x_0}\frac{\mid \mathcal{L}_{n,p}(1+t^2,x)\mid}{(1+x^2)^{1+\alpha}}\\&&+\sup_{x \geq x_0}\frac{\mid f(x) \mid}{(1+x^2)^{1+\alpha}}.
\end{eqnarray*}
The first term of the above inequality tends to zero from Theorem \ref{t2}. By Lemma \ref{l1}(ii), for any fixed $x_0>0$ it is easily seen that $ \sup_{x\geq x_0} \frac{\mid
\mathcal{L}_{n,p}(1+t^2,x)\mid}{(1+x^2)^{1+\alpha}}$
tends to zero as $n \rightarrow \infty$. We can choose $x_0>0$ so
large that the last part of the above inequality can be made small
enough. Thus the proof is completed.
\qed

\section{Voronovskaja type theorem}
In this section we establish a Voronovskaja type asymptotic formula for the operators $\mathcal{L}_{n,p}.$
\begin{lemma} \label{vt1}
For every $x \in [0,\infty)$, we have
\begin{equation}
\lim_{n \rightarrow \infty}n\mathcal{L}_{n,p}((t-x),x)= 1,
\end{equation}
\begin{equation} \label{vt2}
\lim_{n \rightarrow \infty}n\mathcal{L}_{n,p}((t-x)^2,x)= x(2+x).
\end{equation}
\end{lemma}

\begin{theorem}\label{vt3}
If any $f\in C_{x^2}^*[0,\infty)$ such that $f',f''\in C_{x^2}^*[0,\infty)$ and $x \in[0,\infty)$ then, we have
\begin{eqnarray*}
\lim\limits_{n\rightarrow\infty}n\bigg[\mathcal{L}_{n,p}(f,x)-f(x)\bigg] & = & f'(x)+ \big\{(x^2+ 2x)/2\big\}f''(x).\\
\end{eqnarray*}
\end{theorem}
\textbf{Proof.}
Using Taylor's expansion to prove this identity
\begin{eqnarray*}
f(t)-f(x)=(t-x)f'(x)+\frac{1}{2!} f''(x)(t-x)^2+ r(t,x)(t-x)^2,
\end{eqnarray*}
where $r(t,x)$ is Peano form of the remainder, $r(t,x)\in C_{B}[0,\infty)$ and $ \lim_{t
\rightarrow x} r(t,x)=0$. Applying $\mathcal{L}_{n,p}$ to above, we obtain
\begin{equation*}
n[\mathcal{L}_{n,p}(f,x)-f(x)]=f'(x)n\mathcal{L}_{n,p}(t-x,x) +\frac{n}{2!}f''(x)\mathcal{L}_{n,p}((t-x)^2,x)+n \mathcal{L}_{n,p}(r(t,x)(t-x)^2,x).
\end{equation*}
By Cauchy-Schwarz inequality, we have
\newline
\begin{equation} \label{et2}
\mathcal{L}_{n,p}(r(t,x)(t-x)^2,x)\leq
\sqrt{\mathcal{L}_{n,p}(r(t,x)^2,x)}
\sqrt{\mathcal{L}_{n,p}((t-x)^4,x)}.
\end{equation}We observe that $r^2(x,x)=0$ and $r^2(t,x) \in C_{x^2}[0,\infty)$. Then it follows that
\begin{equation} \label{et3}
\lim_{n\rightarrow \infty}n\mathcal{L}_{n,p}(r(t,x)^2,x)=r^2(x,x)=0,
\end{equation}
uniformly with respect to $x\in [0,A]$, where $A>0$. Now from (\ref{et2}), (\ref{et3}) and Lemma \ref{vt1}, we
obtain
\begin{equation*}
\lim_{n\rightarrow \infty}n\mathcal{L}_{n,p}(r(t,x)(t-x)^2,x)=0.
\end{equation*}
Hence,
\begin{eqnarray*}
&&\lim_{n \rightarrow \infty}n[\mathcal{L}_{n,p}(f,x)-f(x)]\\
&& =\lim_{n \rightarrow \infty}\left(f'(x) n\mathcal{L}_{n,p}(t-x,x)+\frac{n}{2}f''(x)\mathcal{L}_{n,p}((t-x)^2,x)
+n \mathcal{L}_{n,p}(r(t,x)(t-x)^2,x)\right)\\
&& = f'(x)+\big\{x(x+2)/2\big\} f''(x),
\end{eqnarray*}
which completes the proof.
\qed

\section{Error Estimation}
The usual modulus of continuity of $f$ on the closed interval $[0, b]$ is defined
by
$$\omega_b(f,\delta) =\sup_{|t-x|\leq\delta,\, x,t\in[0,b]}|f(t)-f(x)|,\,\,  b>0.$$
It is well known that, for a function $f\in E$, $$ \lim_{\delta\rightarrow 0^+}\omega_b(f,\delta)=0,$$
where\\
$$E:=\left\{f\in C[0,\infty):\lim_{x\rightarrow\infty}\frac{f(x)}{1+x^2}\,\, is\,\, finite \right\}.$$
The next theorem gives the rate of convergence of the operators $\mathcal{L}_{n,p}(f,x)$ to  $f(x),$ for all $f \in E.$ 


\begin{theorem}  \label{t2}
Let $f\in E$ and $\omega_{b+1}(f,\delta)$ be its modulus of continuity on the
finite interval $[0,b+1]\subset[0,\infty)$, where $a>0$. Then we have
\begin{equation*}
\| \mathcal{L}_{n,p}(f,x)-f \|_{C[0,b]} \leq
M_f(1+b^2)\delta_n(b)+2\omega_{b+1}\left(f,\sqrt{\delta_n(b)}\right).
\end{equation*}
\end{theorem}
\textbf{Proof.} 
The proof is based on the following inequality
\begin{equation}\label{t3}
\|\mathcal{L}_{n,p}(f,x)-f \| \leq
M_f(1+b^2)\mathcal{L}_{n,p}((t-x)^2,x)+
\left(1+\frac{\mathcal{L}_{n,p}(|t-x|,x)}{\delta}\right)\omega_{b+1}(f,\delta).
\end{equation}
For all $(x,t)\in [0,b]\times[0,\infty):= S.$
To prove (\ref{t3}), we write\\
$$S=S_1\cup S_2:=\{(x,t):0\leq x\leq b,\, 0\leq t \leq b+1\}\cup\{(x,t):0\leq x\leq b,\,  t> b+1\}.$$
If $(x, t)\in S_1,$ we can write
\begin{equation}\label{o}
|f(t)-f(x)|\leq \omega_{b+1}(f,|t-x|)\leq\left(1+\frac{|t-x|}{\delta}\right)\omega_{b+1}(f,\delta)
\end{equation}
where $\delta > 0.$ On the other hand, if $(x, t)\in S_2,$ using the fact that $t-x > 1$,
we have
\begin{eqnarray}\label{o1}
|f(t)-f(x)|&\leq& M_f(1+x^2+t^2)\\
&\leq& M_f(1+3x^2+2(t-x)^2)\nonumber\\
&\leq& N_f(1+b^2)(t-x)^2\nonumber
\end{eqnarray}
where $N_f = 6M_f.$ Combining (\ref{o}) and (\ref{o1}), we get (\ref{t3}).
Now from (\ref{t3}) it follows that
\begin{eqnarray*}
|\mathcal{L}_{n,p}(f,x)-f(x)|&\leq & N_f(1+b^2)\mathcal{L}_{n,p}((t-x)^2,x)+\left(1+\frac{\mathcal{L}_{n,p}(|t-x|,x)}{\delta}\right)\omega_{b+1}(f,\delta)\\
&\leq& N_f(1+b^2)\mathcal{L}_{n,p}((t-x)^2,x)+\left(1+\frac{{[\mathcal{L}_{n,p}((t-x)^2,x)]}^{1/2}}{\delta}\right)\omega_{b+1}(f,\delta).\\
\end{eqnarray*}
By Lemma \ref{l1},  we have
$$\mathcal{L}_{n,p}(t-x)^2\leq\delta_n(b).$$
\begin{equation*}
\| \mathcal{L}_{n,p}(f,x)-f \| \leq
N_f(1+b^2)\delta_n(b)+\left(1+\frac{\sqrt{\delta_n(b)}}{\delta}\right)\omega_{b+1}(f,\delta).
\end{equation*}
Choosing $\delta =\sqrt{\delta_n(b)},$ we get the desired estimation.
\qed

\section{Stancu approach}
In \cite{s} Stancu introduced the following generalization of Bernstein polynomials
\begin{equation}
S_n^\alpha (f,x)=\sum_{k=0}^n f\left(\frac{k}{n}\right)P_{n,\alpha}^k(x),0\leq x \leq 1,
\end{equation}
where $\displaystyle P_{n,\alpha}^k(x)= \binom{n}{k}\frac{\prod_{s=0}^{k-1}(x+\alpha s)\prod_{s=0}^{n-k-1}(1-x+\alpha s)}{\prod_{s=0}^{n-1}(1+\alpha s)}.$\\
\indent We get the classical Bernstein polynomials by putting $\alpha=0$.
Starting with two parameter $\alpha,\beta$ satisfying the
condition $0\leq \alpha \leq \beta$. In 1983, the other
generalization of Stancu opreators was given in \cite{s1} and
studied the linear positive operators $S_n^{\alpha, \beta}:
C[0,1]\rightarrow C[0,1]$ defined for any $f\in C[0,1]$ as
follows:

\begin{equation}
S_n^{\alpha, \beta}(f,x)=\sum_{k=0}^n
p'_{n,k}(x)f\left(\frac{k+\alpha}{n+\beta}\right), 0\leq x \leq 1,
\end{equation}
where $ \displaystyle p'_{n,k}(x)=\binom{n}{k}x^k(1-x)^{n-k}$ is the Bernstein basis
function(cf. \cite{B}).\\

\indent Recently, Ibrahim \cite{i} introduced Stancu-Chlodowsky polynomial and investigated convergence and approximation properties of these operators. Motivated by such types of operators, we introduce a Stancu type generalization of the Baskakov-Schurer-Sz\'{a}sz  type operators (\ref{eq1}) as follows:

\begin{equation}\label{eq2}
\mathcal{L}_{n,p}^{(\alpha,\beta)}(f,x) = (n+p)\sum_{k=0}^\infty b_{n,p}^k(x) \int_0^{\infty} s_{n,p}^{k}(t)f\left(\frac{nt+\alpha}{n+\beta}\right) dt,
\end{equation}
where $b_{n,p}^k(x)$ and $s_{n,p}^{k}(t)$ defined as same in (\ref{eq1}). The operators $\mathcal{L}_{n,p}^{(\alpha,\beta)}(f,x)$ in (\ref{eq2}) are called Baskakov-Schurer-Sz\'{a}sz-Stancu operators. For $\alpha=0, \,\beta=0$ the operators (\ref{eq2}) reduce to the operators (\ref{eq1}).\\
\newline
We know that $$\sum_{k=0}^{\infty} b_{n,p}^{k}(t)=1 ,\text{   } \int_0^{\infty} b_{n,p}^k(t)=\frac{1}{n + p-1},$$  and $$\sum_{k=0}^{\infty} s_{n,p}^k(x)=1 ,\text{   } \int_0^{\infty}
s_{n,p}^k(x)=\frac{1}{n+p}.$$
Before proceeding further, we need some lemmas for proving our next results.

\begin{lemma} \label{sl1} Let $e_m(t) =t^m$, $m = 0,1,2.$ we have
\begin{eqnarray*}
\mathcal{L}_{n,p}^{(\alpha,\beta)}(e_0,x)&&=1,\\
\mathcal{L}_{n,p}^{(\alpha,\beta)}(e_1,x)&& = \frac{(n x+ \alpha)(n+p) + n}{(n+p)(n+\beta)},\\
\mathcal{L}_{n,p}^{(\alpha,\beta)}(e_2,x)&&= \bigg(\frac{n^2 (n+p+1)}{(n+p)(n+\beta)^2}\bigg) x^2 + \left(\frac{4n^2 + 2 n \alpha (n+p)}{(n+p) (n+\beta)^2} \right) x +\bigg[ \frac{2 n^2+ 2 n\alpha (n+p) + \alpha^2 (n+p)^2 }{(n+p)^2(n+\beta)^2}\bigg].
\end{eqnarray*}
\end{lemma}
\textbf{Proof.}
We observe that, $\mathcal{L}_{n,p}^{(\alpha,\beta)}(e_0,x)=\mathcal{L}_{n,p}(e_0,x)=1,$
\begin{eqnarray*}
\mathcal{L}_{n,p}^{(\alpha,\beta)}(e_1,x)&=&\frac{n}{n+\beta}\mathcal{L}_{n,p}(e_1,x)+\frac{\alpha}{n+\beta}\mathcal{L}_{n,p}(e_1,x)\\
&=& \frac{n}{n+\beta}\bigg[x+\frac{1}{n+p}\bigg] +\frac{\alpha}{n+\beta} = \frac{(n x+ \alpha)(n+p) + n}{(n+p)(n+\beta)}.
\end{eqnarray*}
Finally 

\begin{eqnarray*}
\mathcal{L}_{n,p}^{(\alpha,\beta)}(e_2,x)&=&\frac{n^2}{(n+\beta)^2}\mathcal{L}_{n,p}(e_2,x)+
\frac{2n\alpha}{(n+\beta)^2}\mathcal{L}_{n,p}(e_1,x)+\frac{\alpha^2}{(n+\beta)^2}\mathcal{L}_{n,p}(e_0,x)\\
&=&\frac{n^2}{(n+\beta)^2}\left[\bigg(1+\frac{1}{n+p}\bigg) x^2+ \frac{4 x}{(n+p)} +\frac{2}{(n+p)^2}\right]+\frac{2n\alpha}{(n+\beta)^2}\left(x+\frac{1}{n+p}\right)+ \frac{\alpha^2}{(n+\beta)^2}\\
&=& \bigg(\frac{n^2 (n+p+1)}{(n+p)(n+\beta)^2}\bigg) x^2 + \left(\frac{4n^2 + 2 n \alpha (n+p)}{(n+p) (n+\beta)^2} \right) x + \bigg[\frac{2 n^2+ 2 n\alpha (n+p) + \alpha^2 (n+p)^2 }{(n+p)^2(n+\beta)^2}\bigg].
\end{eqnarray*}
\qed

\begin{lemma}
\begin{eqnarray*}
(i) \mathcal{L}_{n,p}^{(\alpha,\beta)}((t-x);x):&=&\mu_{n,p}^1=
\frac{n(1+\alpha)+\alpha p-\beta x (n+p)}{(n+p)(n+\beta)}.\\
(ii) \mathcal{L}_{n,p}^{(\alpha,\beta)}((t-x)^2;x):&=&\mu_{n,p}^2=
\left[\frac{n^2 + \beta^2(n+ p)}{(n+p)(n+\beta)^2}\right] x^2+ \left[\frac{2 n^2 - 2 n \beta (1+\alpha)- 2\alpha p \beta}{(n+p)(n+\beta)^2}\right] x \\&&~~~~+ \left[\frac{ 2 n^2+ 2 n \alpha (n+p) + \alpha^2 (n+p)^2}{(n+\beta)^2 (n+p)^2}\right].
\end{eqnarray*}
\end{lemma}

\begin{theorem}
For each $f \in C_{\rho}^{*}[0,\infty)$, we have
\begin{equation*}\lim _{n\rightarrow \infty} \| \mathcal{L}_{n,p}^{(\alpha,\beta)} (f)-f \|_{\rho} = 0.\end{equation*}
\end{theorem}
\textbf{Proof.}
Using theorem in \cite{AD} we see that it is sufficient to verify the following three conditions
\begin{equation}\label{l6}
\lim_{n\rightarrow \infty} \|
\mathcal{L}_{n,p}^{(\alpha,\beta)}(t^r,x)-x^r\|_{\rho}=0, \text{  } r=0,1,2.
\end{equation}
Since, $\mathcal{L}_{n,p}^{(\alpha,\beta)}(1,x)=1$, the first condition of (\ref{l6}) is satisfied for $r=0$. Now,

\begin{eqnarray*}
\|\mathcal{L}_{n,p}^{(\alpha,\beta)}(t,x)-x\|_{\rho}&=&\sup_{x\in [0,\infty)} \frac{\mid \mathcal{L}_{n,p}^{(\alpha,\beta)}(t,x)-x \mid}{1+x^2}\\
&\leq & \frac{\beta}{(n+\beta)}\sup_{x\in [0,\infty)}\frac{ x}{1+x^2} + \bigg(\frac{n+ \alpha(n+p)}{(n+\beta)(n+p)}\bigg)\\
&\leq & \frac{\beta}{(n+\beta)} + \bigg(\frac{n+ \alpha(n+p)}{(n+\beta)(n+p)}\bigg)
\end{eqnarray*}
\hspace{3.9cm}$\rightarrow 0$ as $n\rightarrow \infty$.\\
Condition of (\ref{l6}) hold for r = 1.
Similarly, we can write

\begin{eqnarray*}
\|\mathcal{L}_{n,p}^{(\alpha,\beta)}(t^2,x)-x^2\|_{\rho} & = & \sup_{x\in [0,\infty)} \frac{\mid \mathcal{L}_{n,p}^{(\alpha,\beta)}(t^2,x)-x^2 \mid}{1+x^2}\\
&\leq & \left[\frac{n^2(n+p+1)}{(n+p)(n+\beta)^2}-1 \right] \sup_{x\in [0,\infty)} \frac{x^2}{1+x^2}+  \left(\frac{4n^2 + 2 n \alpha (n+p)}{(n+p) (n+\beta)^2} \right) \sup_{x\in [0,\infty)} \frac{x}{1+x^2} \\&&+ \bigg[ \frac{2 n^2+ 2 n\alpha (n+p) + \alpha^2 (n+p)^2 }{(n+p)^2(n+\beta)^2}\bigg].
\end{eqnarray*}
which implies that \\
$$\|\mathcal{L}_{n,p}^{(\alpha,\beta)} (t^2;x)- x^2 \|_{\rho} = 0$$ as $n\rightarrow \infty $. Thus, from Gadzhiev’s Theorem in \cite{AD} we obtain the desired result of theorem. 
\qed

\begin{theorem}  
Let $f \in C_{\rho}[0,\infty)$ and $\omega_{[0,a+1]}(f,\delta)$ be its modulus of continuity on the
finite interval $[0,a+1]\subset[0,\infty)$, where $a>0$. Then, we have
\begin{eqnarray*}
\| \mathcal{L}_{n,p}^{(\alpha,\beta)}(f)-f \|_{[0,a]} \leq C\bigg[(1+a^2)\delta_n(x)+2~ \omega_{[0,a+1]}\left(f,\sqrt{\delta_n(x)}\right)\bigg], ~~~ C > 0.
\end{eqnarray*}
\end{theorem}
\textbf{Proof.} The proof is similar to Theorem \ref{t2} with appropriate modifications. We
omit the details. \qed

\section{Rate of convergence}

We can give some estimations of the errors
$|\mathcal{L}_{n,p}^{(\alpha,\beta)}(f)-f |,$ $n\in \mathbb{N}$ for unbounded functions
by using a weighted modulus of smoothness associated to the space $B_{\rho_{\gamma}}{(\mathbb{R}_+)}$.
The weighed modulus of continuity $\Omega_{\rho_{\gamma}}(f;\delta)$ was defined by L\'{o}pez--Moreno in \cite{LM}.
We consider
\begin{equation}\label{r1}
\Omega_{\rho_{\gamma}}(f;\delta) =\sup_{x\geq 0, 0\leq h\leq \delta} \frac{|f(x+h)-f(x)|}{1+(x+h)^{2+\gamma}},\,\,  \delta>0, \gamma\geq 0.
\end{equation}
It is evident that for each $f\in B_{\rho_{\gamma}}{(\mathbb{R}_+)}, \  \Omega_{\rho_{\gamma}}(f; \cdot)$ is well defined and

$$\Omega_{\rho_{\gamma}}(f;\delta)\leq 2\| f \|_{\rho_{\gamma}}.$$
The weighted modulus of smoothness $\Omega_{\rho_{\gamma}}(f; \cdot)$ possesses
the following properties.
\begin{align} \label{p1}
&(i) \Omega_{\rho_{\gamma}}(f; \lambda \delta)\leq (\lambda+1) \Omega_{\rho_{\gamma}}(f;\delta),\delta>0, \lambda>0 \\
&(ii) \Omega_{\rho_{\gamma}}(f; n \delta)\leq n\Omega_{\rho_{\gamma}}(f; \delta), ~~n\in\mathbb{N}\nonumber\\
&(iii) \lim\limits_{\delta \rightarrow  0}\Omega_{\rho_{\gamma}}(f; \delta)= 0.\nonumber
\end{align}
Now, we are ready to prove our next theorem by using above properties.  

\begin{theorem}
For all non-decreasing $f \in B_{\rho_{\gamma}}{(\mathbb{R}_+)}$, we have
\begin{equation*}
| \mathcal{L}_{n,p}^{\alpha,\beta} (f,x)-f| \leq \sqrt{\mathcal{L}_{n,p}^{(\alpha,\beta)}(\nu^2_{x,\gamma};x)}\left( 1 + \frac{1}{\delta}\sqrt{\mathcal{L}_{n,p}^{(\alpha,\beta)}(\Psi_{x}^2;x)}\right)\Omega_{\rho_{\gamma}}(f; \delta),
\end{equation*} 
$x\geq 0,~ \delta> 0, ~ n\in \mathbb{N},$ where $$ \nu_{x,\gamma}(t):=1+(x+|t-x|)^{2+\gamma}, ~~ \Psi_{x}(t):=|t-x|, ~ t\geq 0.$$
\end{theorem}
\textbf{Proof.}
Let $n\in\mathbb{N}$ and $f \in B_{\rho_\gamma }(\mathbb{R_+}).$ From (\ref{r1}) and (\ref{p1}), we can write

\begin{eqnarray*}
|f(t)-f(x)| &\leq& \bigg(1+(x + |t-x|)^{2+\gamma}\bigg)\bigg(1+\frac{1}{\delta} |t-x|\bigg)\Omega_{\rho_{\gamma}}(f; \delta)\\
&=&\nu_{x,\gamma}(t)\bigg(1+\frac{1}{\delta}\Psi_x(t)\bigg)\Omega_{\rho_{\gamma}}(f; \delta).
\end{eqnarray*}
Now, applying operator $\mathcal{L}_{n,p}^{(\alpha,\beta)}$ on above inequality, we get

\begin{eqnarray*}
| \mathcal{L}_{n,p}^{(\alpha,\beta)} (f,x)-f(x)|
&\leq & \Omega_{\rho_{\gamma}}(f; \delta)\mathcal{L}_{n,p}^{(\alpha,\beta)}\bigg(\nu_{x,\gamma}\bigg(1+\frac{1}{\delta}\Psi_x\bigg);x\bigg)
\end{eqnarray*}
\begin{eqnarray} \label{p2}
~~~~~~~~~~~~~~~~~~~~~~~~~\hspace{2cm}&\leq &\Omega_{\rho_{\gamma}}(f; \delta) \bigg(\mathcal{L}_{n,p}^{(\alpha,\beta)}(\nu_{x,\gamma};x)+\mathcal{L}_{n,p}^{(\alpha,\beta)}\bigg(\frac{\nu_{x,\gamma} \Psi_x}{\delta};x\bigg)\bigg).
\end{eqnarray}
By using the Cauchy-Schwartz inequality, we obtain
\begin{eqnarray*}
\mathcal{L}_{n,p}^{(\alpha,\beta)}\bigg(\frac{\nu_{x,\gamma} \Psi_{x}}{\delta};x\bigg) &\leq & \bigg\{\mathcal{L}_{n,p}^{(\alpha,\beta)}\big(({\nu_{x,\gamma})^2;x)}\bigg\}^{1/2} \bigg\{\mathcal{L}_{n,p}^{(\alpha,\beta)}\bigg(\bigg(\frac{\Psi_x}{\delta}\bigg)^2;x\bigg)\bigg\}^{1/2}\\
&=&\frac{1}{\delta}\bigg\{\mathcal{L}_{n,p}^{(\alpha,\beta)}({\nu^2_{x,\gamma};x)}\bigg\}^{1/2} \bigg\{\mathcal{L}_{n,p}^{(\alpha,\beta)}({\Psi_x}^2;x)\bigg\}^{1/2}.
\end{eqnarray*}
Now, by (\ref{p2}), we get
\begin{equation*}
| \mathcal{L}_{n,p}^{(\alpha,\beta)} (f,x)-f| \leq \sqrt{\mathcal{L}_{n,p}^{(\alpha,\beta)}(\nu^2_{x,\gamma};x)}\left( 1 + \frac{1}{\delta}\sqrt{\mathcal{L}_{n,p}^{(\alpha,\beta)}(\Psi_{x}^2;x)}\right)\Omega_{\rho_{\gamma}}(f; \delta).
\end{equation*} 

\section{Better approximation}
It is well know that the operators preserve constant as well as linear functions. To make the convergence faster, King \cite{kin} proposed an approach  to modify the classical Bernstein polynomials, so that this sequence preserves two test functions $e_0$ and $e_1$.
 After this several researchers have studied that many approximating operators $L$, possess these properties i.e.
 $ L(e_i,x)=e_i(x)$ where $e_i(x)=x^i(i=0,1)$, for examples Bernstein, Baskakov and Baskakov-Durrmeyer-Stancu operators.\\
\indent As the operators $\mathcal{L}_{n,p}^{(\alpha,\beta)}$  introduced in (\ref{eq2}) preserve only the constant functions so further modification of these operators is proposed to be made so that the modified operators preserve the constant as well as linear functions, for this purpose the modification of $\mathcal{L}_{n,p}^{(\alpha,\beta)}$ as follows:

\begin{equation*}
\mathcal{L}_{n,p}^{*(\alpha,\beta)}(f,x) = (n+p)\sum_{k=0}^\infty b_{n,p}^k (r_n(x)) \int_0^{\infty/A} s_{n,p}^{k}(t)f\left(\frac{n t+\alpha}{n+\beta}\right) dt
\end{equation*}
where $r_n(x)=\bigg(\frac{(n+\beta)x-\alpha}{n}-\frac{1}{n+p}\bigg)$ and $x\in I_n=\bigg[\frac{n+\alpha(n+p)}{(n+\beta)(n+p)},\infty \bigg)$.
\begin{lemma}
For each $x\in I_n$, we have
$$\mathcal{L^{*}}_{n,p}^{(\alpha,\beta)}(1,x)=1 ,\text{ } \mathcal{L^{*}}_{n,p}^{(\alpha,\beta)}(t,x) = x,$$
\begin{eqnarray*}
\mathcal{L^{*}}_{n,p}^{(\alpha,\beta)}(t^2,x)&=& \bigg(\frac{ n^2(2(n+p)+1)+ (n+p)(n+\beta)^2}{(n+p)(n+\beta)^2} \bigg) x^2 + \left(\frac{6 n^2 + 4 n \alpha (n+p)}{(n+p) (n + \beta)^2} \right)x \\ &&+
\left(\frac{3 n^2 + 4 n \alpha (n+p) + 2 \alpha^2(n+p)^2}{(n+p)^2 (n+\beta)^2}\right).
\end{eqnarray*}
\end{lemma}

\begin{lemma}
For $x\in I_n$, the following holds,
$$\widetilde{\mu}_{n,p}^{1}(x)=\mathcal{L^{*}}_{n,p}^{(\alpha,\beta)}(t-x,x)=0,$$
\begin{eqnarray*}
\widetilde{\mu}_{n,p}^{2}(x) = \mathcal{L^{*}}_{n,p}^{(\alpha,\beta)}((t-x)^2,x)&=& \bigg(\frac{ n^2(2(n+p)+1)}{(n+p)(n+\beta)^2} \bigg) x^2 + \left(\frac{6 n^2 + 4 n \alpha (n+p)}{(n+p) (n + \beta)^2} \right)x \\ &&+
\left(\frac{3 n^2 + 4 n \alpha (n+p) + 2 \alpha^2(n+p)^2}{(n+p)^2 (n+\beta)^2}\right).
\end{eqnarray*}
\end{lemma}

\begin{theorem}
Let $f\in C_B(I_n)$ and $x\in I_n$. Then, there exist an absolute constant $C>0$ such that
$$\bigg|  \mathcal{L^{*}}_{n,p}^{(\alpha,\beta)}(f,x)-f(x)\bigg| \leq C\omega_2\left(f,\sqrt{\widetilde{\mu}_{n,p}^{2}(x)} ~\right).$$
\end{theorem}
\textbf{Proof.}
Let $g \in C_B(I_n)$ and $x,t \in I_n$. By Taylor's expansion we have
\begin{equation}
g(t)=g(x)+ (t-x)g'(x)+\int_x^t(t-u)g''(u)du.
\end{equation}
Applying $\mathcal{L^{*}}_{n,p}^{(\alpha,\beta)}$, we get
\begin{equation*}
\mathcal{L^{*}}_{n,p}^{(\alpha,\beta)}(g,x)-g(x) = g'(x)\mathcal{L^{*}}_{n,p}^{(\alpha,\beta)}((t-x),x)+\mathcal{L^{*}}_{n,p}^{(\alpha,\beta)}\left(\int_x^t(t-u)g''(u)du,x\right).
\end{equation*}
Now, we have $\displaystyle \bigg|\int_x^t(t-x)g''(u)du\bigg| \leq (t-x)^2\|g''\|$,
$$ \bigg| \mathcal{L^{*}}_{n,p}^{(\alpha,\beta)}(g,x)-g(x)\bigg| \leq \mathcal{L^{*}}_{n,p}^{(\alpha,\beta)}((t-x)^2,x)\|g''\|=\widetilde{\mu}_{n,p}^{2}\|g''\|.$$
Since $\bigg|\mathcal{L^{*}}_{n,p}^{(\alpha,\beta)}(f,x)\bigg|\leq \|f\|$,
\begin{eqnarray*}
\bigg| \mathcal{L^{*}}_{n,p}^{(\alpha,\beta)}(f,x)-f(x)\bigg| &\leq& \bigg| \mathcal{L^{*}}_{n,p}^{(\alpha,\beta)} (f-g,x)-(f-g)(x)\bigg| +\bigg|\mathcal{L^{*}}_{n,p}^{(\alpha,\beta)}(g,x)-g(x)\bigg|\\
&\leq& 2\|f-g\|+\widetilde{\mu}_{n,p}^{2}\|g''\|.
\end{eqnarray*}
Taking infimum overall $g\in C^2(I_n)$, we obtain
$$\bigg|\mathcal{L^{*}}_{n,p}^{(\alpha,\beta)}(f,x)-f(x)\bigg|\leq K_2(f,\widetilde{\mu}_{n,p}^{2}).$$
In view of (\ref{1}), we have
$$\bigg|\mathcal{L^{*}}_{n,p}^{(\alpha,\beta)}(f,x)-f(x)\bigg| \leq C \omega_2\left(f,\sqrt{\widetilde{\mu}_{n,p}^{2}}\right),$$
which proves the theorem.


\section{$q$-Baskakov-Schurer-Sz\'{a}sz operator}
Let $p, k, n\in\mathbb{N}, A > 0$ and $f$ be a real valued continuous function on the interval $[0,\infty)$. In \cite{IY} Y\"{u}ksel defined the $q$-Baskakov-Schurer-Sz\'{a}sz type linear positive operators as 

\begin{equation}\label{qb}
\mathcal{L}^{q}_{n,p}(f,x)= [n+p]_q \sum_{k=0}^{\infty} b^k_{n,p}(x;q)\int_{0}^{\frac{\infty}{A(1-q)}} f(t) s^k_{n,p}(t;q) d_q t 
\end{equation}
where 
$$\displaystyle b^k_{n,p}(x;q)=\left[\begin{array}{c} n+p+k-1\\k \end{array}\right]_q q^{k^2}\frac{x^k}{(1+x)_q^{n+k+p}},$$
and\,\, $\displaystyle s^k_{n,p}(t;q)= e^{-[n+p]_qt}\frac{{([n+p]_qt)}^k}{k!}.$\\
\newline
If $q = 1$ then the operators $\mathcal{L}^q_{n,p}(f,x)$
are reduced to Baskakov-Schurer-Sz\'{a}sz type operators defined in (\ref{eq1}).

\begin{lemma}\cite{IY}\label{ql1}
Let $e_m(t) =t^m$, $m = 0,1,2.$ we have
\begin{itemize}
\item[(i)] $\mathcal{L}^q_{n,p}(e_0;x)=1,$
\item[(ii)] $\mathcal{L}^q_{n,p}(e_1;x)=\frac{x}{q^2}+\frac{1}{q[n+p]_q},$
\item[(iii)] $\mathcal{L}^q_{n,p}(e_2;x)=\frac{[n+p+1]_q}{q^{6}[n+p]_q}x^2+\frac{1+2q+q^2}{q^5[n+p]_q}x+\frac{1+q}{q^3[n+p]^2_q}.$
\end{itemize}
\end{lemma}

\begin{lemma}\cite{IY}\label{al1}
Let $ q\in (0,1)$ and  $p\in \mathbb{N}$. Then following hold:\\
$\delta_{n,p}(x;q):=\mathcal{L}^q_{n,p}((t-x)^2;x)\leq \frac{9}{q^6}\left(1-q^3+\frac{1}{[n+p]_q}\right)(x+1)^2.$
\end{lemma}

\begin{remark}\label{rem1}
$\alpha_{n,p}(x;q):=\mathcal{L}^q_{n,p}((t-x);x)=\bigg(1-\frac{1}{q^2}\bigg) x+\frac{1}{q[n+p]_q}.$
\end{remark}
Now, we establish a Voronovskaja type asymptotic formula for the operators $\mathcal{L}^q_{n,p}.$

\begin{theorem}\label{vothm2}
Let $f$  be bounded and integrable on the interval $[0,\infty).$ First and second derivatives of $f$ exists at a fixed point $x\in [0,\infty)$  and $q=q_n\in (0,1)$ such that $q=q_n \rightarrow 1$ as $n\rightarrow \infty$, then
\begin{eqnarray*}
\lim\limits_{n\rightarrow\infty} [n]_q\big[\mathcal{L}^q_{n,p}(f,x)-f(x)\big] & = & f'(x) +  x f''(x).
\end{eqnarray*}
\end{theorem}
\textbf{Proof.}
In order to prove this identity we use Taylor's expansion
\begin{eqnarray*}
f(t)-f(x)=(t-x)f'(x) + (t-x)^2\bigg( \frac{1}{2!} f''(x) + \xi(t-x)\bigg),
\end{eqnarray*}
where $\xi$ is bounded and  $ \lim_{t\rightarrow 0} \xi(t)=0$. Applying $\mathcal{L}^q_{n,p}$ to above, we obtain
\begin{eqnarray*}
[n]_q[\mathcal{L}^q_{n,p}(f,x)-f(x)]&=& f'(x)[n]_q\mathcal{L}^q_{n,p}(t-x,x) +\frac{[n]_q}{2!}f''(x)\mathcal{L}^q_{n,p}((t-x)^2,x)+ [n]_q \mathcal{L}^q_{n,p}(\xi(t-x)(t-x)^2,x)\\
&&= [n]_q f'(x)\alpha_n(q_n,x)+ [n]_q\frac{1}{2} f''(x)\delta_n(q_n,x)+ [n]_q \mathcal{L}^q_{n,p}(\xi(t-x)(t-x)^2,x)
\end{eqnarray*}
where $\alpha_n(q_n,x)$ and $\delta_n (q, x)$ defined as in remark \ref{rem1} and  lemma \ref{al1}.
Using Cauchy-Schwarz inequality, we have

\begin{equation} \label{e2}
[n]_q \mathcal{L}^q_{n,p}(\xi(t-x)(t-x)^2,x)\leq
\sqrt{\mathcal{L}^q_{n,p}(\xi(t-x)^2,x)}
\sqrt{ [n]_q^2 \mathcal{L}^q_{n,p}((t-x)^4,x)}.
\end{equation}
We observe that $\lim_{t\rightarrow x}\xi^2(t,x)=0$ it follows that $\lim_{[n]_q\rightarrow \infty}[n]_q\mathcal{L}^q_{n,p}(\xi(t-x)(t-x)^2,x)=0$


Hence,
\begin{eqnarray*}
&&\lim_{[n]_q \rightarrow \infty} [n]_q [\mathcal{L}^q_{n,p}(f,x)-f(x)]\\
&& =\lim_{n \rightarrow \infty}\left(f'(x) [n]_q\mathcal{L}^q_{n,p}(t-x,x)+\frac{[n]_q}{2}f''(x)\mathcal{L}^q_{n,p}((t-x)^2,x)
+[n]_q \mathcal{L}^q_{n,p}(\xi(t-x)(t-x)^2,x)\right)\\
&& = f'(x)+ x f''(x),
\end{eqnarray*}
which completes the proof.
\qed

\begin{theorem}
For $q_n\in (0, 1),$ the sequence $\{\mathcal{L}^q_{n,p}\}$
 converges to $f$ uniformly on $[0,A],~ A>0$ for each $f \in C^{*}_{x^2} [0,\infty)$ if and only if $\lim_{n\rightarrow\infty} q_n = 1.$
\end{theorem}
\textbf{Proof.}  The proof is similar to the Theorem 2 in \cite{VG}.

\section{Weighted statistical approximation}
Our next concern is the study of statistical convergence of the sequence of the $q$-Baskakov-Schurer-Sz\'{a}sz operators. 
Now, we recall the concept of statistical convergence for
sequences of real numbers which was introduced by Fast \cite{HF} and further studied by many others.

Let $K\subseteq \mathbb{N}$ and $K_{n}=\left\{ j\leq n:j\in
K\right\}.$ Then the $natural~density$ of $K$ is defined by $\delta
(K)={\lim\limits_{n}}~ n^{-1}|K_{n}|$ if the limit exists, where $|K_{n}|$ denotes the
cardinality of the set $K_{n}$.

\parindent=8mm A sequence $x=(x_{j})_{j\geq1}$ of real numbers is said to be $%
statistically$ $convergent$ to $L$ provided that for every $\epsilon >0$ the
set $\{j\in \mathbb{N}:|x_{j}-L|\geq \epsilon \}$ has natural density zero,
i.e. for each $\epsilon >0$,
\begin{equation*}
\lim\limits_{n}\frac{1}{n}|\{j\leq n:|x_{j}-L|\geq \epsilon \}|=0.
\end{equation*}
It is denoted by $st-\lim\limits_{n}x_{n}=L$.\\
\newline
We consider a sequence $q=(q_{n}),$ $q_{n}\in $ $(0,1),$ such that
\begin{equation}\label{s1}
\lim\limits_{n\rightarrow \infty }q_{n}=1.
\end{equation}
The condition (\ref{s1}) guarantees that $[n]_{q_{n}}\rightarrow
\infty $ as $n\rightarrow \infty .$\\
\indent In \cite{DO} Doğru and Kanat defined the Kantorovich-type modification of Lupa\c{s} operators as follows:
\begin{equation}\label{dk}
\tilde{R}_n(f;q;x)=[n+1]\sum_{k=0}^{n}\bigg(\int_{\frac{[k]}{[n+1]}}^{\frac{[k+1]}{[n+1]}} ~ f(t) d_{q}t\bigg)\left(\begin{array}{c}n \\k \end{array}\right) \frac{q^{-k}q^{k(k-1)/2} x^k (1-x)^{(n-k)}}{(1-x+qx)\cdots(1-x+q^{n-1} x)}.
\end{equation}
Doğru and Kanat \cite{DO} proved the following statistical Korovkin-type approximation
theorem for operators (\ref{dk}).

\begin{theorem}
Let $q:=(q_n),~ 0 < q < 1$, be a sequence satisfying the following conditions:\\
\begin{equation}\label{thm10}
 st-\lim_n q_n = 1,~ st - \lim_n q_n^n = a ~( a < 1)~ and~ st - \lim_n \frac{1}{[n]_q} = 0,
\end{equation}
then if $f$ is any monotone increasing function defined on $[0, 1]$, for the positive linear operator $\tilde{R}_n(f;q;x)$, then
$$ st- \lim_n {\Vert \tilde{R}_n(f;q;\cdot) - f \Vert}_{C[0, 1]} = 0$$ holds. 
\end{theorem}
In \cite{do26} Do\u{g}ru gave some examples so that $(q_{n})$ is statistically convergent to $1$ but it may not convergent to $1$ in the ordinary case.

\begin{theorem} 
Let ${{\mathcal{L}^{q}_{n,p}}}$ be the sequence of the operators (\ref{qb}) and the sequence
$q = (q_n)$ satisfies (\ref{thm10}). Then for any function $f \in C[0,\nu] \subset C[0,\infty),~ \nu > 0 ,$ we have
\begin{equation}
st- \lim_n \Vert \mathcal{L}^q_{n,p}(f; \cdot) - f \Vert = 0,
\end{equation}
where $C[0,\nu]$ denotes the space of all real bounded functions $f$ which are continuous in
$[0,\nu].$
\end{theorem}
\textbf{Proof.} The proof is similar to Theorem 2 \cite{KJ} with appropriate modifications. So, we
omit the details.\\

\begin{theorem}\label{thm3}
Let $\mathcal{L}^q_{n,p}$ be the sequence of the operators (\ref{qb}) and the sequence $q=(q_{n})$ satisfies (\ref{thm10}). Then for any function $f\in C_{B}[0,\infty ),$ we have

\begin{equation*}
st-\lim_{n\rightarrow \infty }{\Vert }\mathcal{L}^q_{n,p}(f;.)-f{\Vert }%
_{\rho}=0.
\end{equation*}
\end{theorem}

\hspace{-0.9cm}
\textbf{Proof.} By Lemma (\ref{ql1})(iii), we have  $\mathcal{L}^q_{n,p}(t^2; x) \leq Cx^2,$ where $C$ is a positive constant, $\mathcal{L}^q_{n,p}(f; x)$ is a sequence of positive linear operator acting from $C_{\rho}[0,\infty)$ to $B_{\rho}[0,\infty)$. \\
Using $\mathcal{L}^q_{n,p}(1;x) = 1,$ it is clear that\\

$$st-\lim_{n}\|\mathcal{L}^q_{n,p}(1;x)-1\|_{\rho}=0.$$
Now, by Lemma (\ref{ql1})(ii), we have

$$\|\mathcal{L}^q_{n,p}(t; x)-x\|_{\rho}= \sup_{x\in [0,\infty)} \frac{|\mathcal{L}^q_{n,p}(t;x)-x|}{1+x^2} \leq \bigg( \frac{1-q^2}{q^2}\bigg) \sup_{x\in [0,\infty)} \frac{x}{1+x^2} + \frac{1}{q [n+p]_q}.$$
Using (\ref{thm10}), we get

$$st-\lim_{n}\left( \bigg( \frac{1-q^2}{q^2}\bigg)  + \frac{1}{q [n+p]_q }\right)=0,$$
then $$st-\lim_{n}\|\mathcal{L}^q_{n,p}(t;x)-x\|_{\rho}=0.$$
Finally, by Lemma(\ref{ql1})(iii), we have

\begin{align*}
\Vert\mathcal{L}^q_{n,p}(f,x)(t^2;q_{n};x)-x^2\Vert_{\rho}
&\leq \left({\frac{[n+p+1]_q}{q^6 [n+p]_q}-1} \right)\sup_{x\in[0,\infty)}\frac{x^2}{1+x^2}\\
&~~~+\left(\frac{1+2q+q^2}{q^5[n+p]_q} \right)\sup_{x\in[0,\infty)}\frac{x}{1+x^2} + \frac{1+q}{q^3[n+p]^2_q}.\\
& \leq \left({\frac{[n+p+1]_q}{q^6 [n+p]_q}-1} \right) +\left(\frac{1+2q+q^2}{q^5[n+p]_q} \right)+ \frac{1+q}{q^3[n+p]^2_q}.\\
\end{align*}
Now, If we choose

\begin{equation*}
\alpha _{n}=  \left({\frac{[n+p+1]_q}{q^6 [n+p]_q}-1} \right),
\end{equation*}%
\begin{equation*}
\beta _{n}=\left(\frac{1+2q+q^2}{q^5[n+p]_q} \right),
\end{equation*}%
\begin{equation*}
\gamma _{n}= \frac{1+q}{q^3[n+p]^2_q},
\end{equation*}%
then by $(\ref{thm10}),$ we can write

\begin{equation}\label{3.1}
st-\lim_{n\rightarrow \infty }\alpha _{n}=0=st-\lim_{n\rightarrow \infty
}\beta _{n}=st-\lim_{n\rightarrow \infty }\gamma _{n}.
\end{equation}
Now for given $\epsilon >0$, we define the following four sets:

\begin{equation*}
{S}=\{k: \Vert \mathcal{L}^q_{n,p}(t^{2};q_{k};x)-x^{2}\Vert_{\rho} \geq \epsilon \},
\end{equation*}%
\begin{equation*}
S_{1}=\lbrace k:\alpha _{k}\geq \frac{\epsilon}{3}\rbrace,
\end{equation*}%
\begin{equation*}
S_{2}=\{k:\beta _{k}\geq \frac{\epsilon }{3}\},
\end{equation*}%
\begin{equation*}
S_{3}=\{k:\gamma _{k}\geq \frac{\epsilon }{3}\}.
\end{equation*}%
It is obvious that $S  \subseteq {S_{1}\cup S_{2}\cup S_{3}}$. Then, we obtain%
\begin{eqnarray*}
\delta \big(\{k &\leq &n:\Vert \mathcal{L}^q_{n,p}(t^{2};x)-x^{2}\Vert_{\rho} \geq
\epsilon \}\big) \\
&&~~~~~~\leq \delta\big(\{k\leq n:\alpha _{k}\geq \frac{\epsilon }{3}\}\big)+\delta
\big(\{k \leq n:\beta _{k}\geq \frac{\epsilon }{3}\}\big)+\delta\big(\{k\leq n:\gamma
_{k}\geq \frac{\epsilon }{3}\}\big).
\end{eqnarray*}%
\newline
Using (\ref{3.1}), we get%
\begin{equation*}
st-\lim_{n\rightarrow \infty }\Vert \mathcal{L}^q_{n,p}(t^{2};x)-x^{2}\Vert_{\rho} =0.
\end{equation*}
Since
\begin{eqnarray*}
\Vert \mathcal{L}^q_{n,p}(f;x)-f\Vert_{\rho} \leq \Vert \mathcal{L}^q_{n,p}(t^{2};x)-x^{2}\Vert_{\rho} +\Vert  \mathcal{L}^q_{n,p}(t;x)-x\Vert_{\rho}+{\Vert \mathcal{L}^q_{n,p}(1;x)-1\Vert}_{\rho} ,\newline
\end{eqnarray*}
we get
\begin{eqnarray*}
st-\lim_{n\rightarrow \infty }\Vert\mathcal{L}^q_{n,p}(f;x)-f\Vert_{\rho} 
&\leq &st-\lim_{n\rightarrow \infty }\Vert\mathcal{L}^q_{n,p}(t^{2};x)-x^{2}\Vert_{\rho}\\&& +~st-\lim_{n\rightarrow \infty }\Vert \mathcal{L}^q_{n,p}(t;x)-x\Vert_{\rho}\\
&& + ~st-\lim_{n\rightarrow \infty }\Vert \mathcal{L}^q_{n,p}(1;x)-1\Vert_{\rho} ,
\end{eqnarray*}%
which implies that
\begin{equation*}
st-\lim_{n\rightarrow \infty }\Vert \mathcal{L}^q_{n,p}(f;x)-f\Vert_{\rho} =0.
\end{equation*}
This completes the proof of the theorem.\qed


\vspace{0.6cm}
\hspace{-0.9cm}\textbf{Acknowledgment}\\
The authors are thankful to the referees for making valuable suggestions.

\vspace{0.9cm}
\hspace{-0.9cm}

\end{document}